\documentclass[12pt, a4paper]{article}

\usepackage{latexsym}
\usepackage{a4wide}
\usepackage{amscd}
\usepackage{graphics}
\usepackage{amsmath}
\usepackage{amssymb}
\usepackage[latin1]{inputenc}
\usepackage{mathrsfs}
\input xy
\xyoption{all}
\usepackage{bm}

\newcommand{\C}{\mathbb C}

\newcommand{\R}{\mathbb R}
\newcommand{\RR}{{\mathbb R}}

\newcommand{\Z}{\mathbb Z}
\newcommand{\sgn}{\mathrm{sign}}
\newcommand{\Ind}{\mathrm{Ind\,}}
\newcommand{\spfl}{\mathrm{sf}}
\newcommand{\dd}{\mathrm d}
\newcommand{\Id}{\mathrm{Id}}
\newcommand{\Imm}{\mathrm{Im\,}}
\newcommand{\Gl}{\mathrm{Gl}}
\newcommand{\dirlim}{\mathrm{dirlim}}
\newcommand{\Bif}{\mathrm{Bif\,}}

\newcommand{\KK}{K}

\newcommand{\Vect}{{\mathrm{Vect}}}
\newcommand{\Fred}{{\mathscr F}}

\newcommand{\TF}{{\mathrm{TF}}}
\newcommand{\finedim}{\hfill $q.e.d.$\\}            

\newcommand{\proof}{{\sl Proof.}\hspace{5pt}}   

\newcommand{\Ddx}{\tfrac{ D}{dx}}
\newcommand{\ga}{\gamma}
\newcommand{\ispec}{\mu_{\scriptscriptstyle{\mathrm{spec}}}}
\newcommand{\Ddtt}{\tfrac{ D^2}{ dt^2}}
\newcommand{\la}{\lambda}
\newcommand{\icon}{\mu_{\scriptscriptstyle{\mathrm{con}}}}

\newtheorem{mainthm}{\sc Theorem}           
\newtheorem{thm}{\sc Theorem}[section]      
\newtheorem{prop}[thm]{\sc  Proposition}     
\newtheorem{defn}[thm]{\sc Definition}      
\newtheorem{rem}[thm]{\sc Remark}       
\newtheorem{notation}{\sc Notation}    
\begin{document}

\title{A $K$-theoretical invariant and bifurcation for a parameterized family of
functionals}

\author{Alessandro Portaluri \thanks{ The author was partially
supported by PRIN {\em Variational Methods and Nonlinear
Differential Equations.}}\\ \small{Dipartimento di Matematica {\em
Ennio De Giorgi\/},}\\ \small{Ex-Collegio Fiorini,  Università del Salento,}\\
\small{Via per Arnesano, Lecce, Le, Italy.}\\\small{ Email:
alessandro.portaluri@unile.it.}}

\date{\today}
\maketitle
\tableofcontents




\begin{abstract}
Let $\mathscr F:= \{f_x:x \in X\}$ be a  family of functionals
defined on a Hilbert manifold $\widetilde E$ and smoothly
parameterized by  a compact connected orientable  $n$-dimensional
manifold $X$, and let $\sigma:X \to \widetilde E$ be a smooth
section of critical points of $\mathscr F$. The aim of this paper is
to give a sufficient topological condition on the parameter space
$X$ which detects bifurcation of critical points for $\mathscr F$
from the trivial branch. Finally we are able to give some
quantitative properties of the  bifurcation set for perturbed
geodesics on semi-Riemannian manifolds.
\end{abstract}

\section{Introduction}

The meaning of the word {\em bifurcation\/} changes depending on the
context and here we deal only with bifurcation from the {\em trivial
branch\/}. This is perhaps the simplest bifurcation phenomena which
occur in nonlinear analysis, differential geometry and mechanics.
The aim of this paper is to establish an abstract bifurcation result
in a form suitable for the applications to geometric variational
problems.

Given two normed linear spaces $E$ and $F$ and a continuous family
of Fredholm maps $h\colon X \times E \to F$ continuously
parameterized by a path connected topological space such that
$h(\cdot,0)=0$, we say that a point $x_* \in X$ is a {\em
bifurcation point from the trivial branch $X \times \{0\}$\/} if
each neighborhood of $(x_*,0)$ contains solutions of $h(x,v)=0$ with
$v \not=0$. It is a well-known fact that there exists a strict
relation between the topology of the parameter space and the
existence of a bifurcation point from the trivial branch and thus an
important role is played by the homotopy theoretical methods. The
literature on this topic is quite broad and the main contributions
were given  in the past decades mainly  by Alexander, Ize,
Fitzpatrick, Pejsachowicz among others. (See for instance
\cite{FitPej91,Pej88,Pej01} and references therein). The Alexander's
approach can be quickly sketched as follows. Consider a family $f_p
: \R^n \to \R^n$ parameterized by $p \in \R^k$, with $f_p(0)=0$ and
let $L_p$ be the derivative of $f_p$ at $0$. Let $p_0$ be an
isolated point in the set $\Sigma(f)$ of all parameters at which the
derivative $L_p$ is singular. Restricting the map $p \mapsto L_p$ to
the boundary of a small disk $D$ centered at $p_0$ one gets an
element $\gamma_f$ in the homotopy group $\pi_{k-1}(\Gl(n)).$ This
element is an obstruction to the existence of a deformation of the
map $L:(D,
\partial D)\to (M_n(\R), \Gl(n))$ to a map into $\Gl(n)$. The point
$p_0$ is a bifurcation point, provided the image of $\gamma_f$ by
the classical Whitehead $J$-homomorphism is nontrivial.

However the inspiration for our results is a series of papers
authored by Fitzpatrick and Pejsachowicz \cite{FitPej91},
\cite{Pej88} \cite{Pej01} among others. In this paper the authors
defined a bifurcation invariant which could be roughly seen as a
global version of the Alexander's invariant. In fact, given a family
of quasilinear Fredholm maps of index zero continuously
parameterized by a smooth compact manifold $X$ with $f_x(0)=0$, the
obstruction to deforming such a family to a family of isomorphism is
the {\em analytical index\/} (or {\em index bundle\/}), namely an
element of the  Grothendieck group of virtual vector bundles
$\KK^0(X)$ and quite naturally the bifurcation invariant in this
setting is provided by $J(\Ind L)$, where $J$ is the
$J$-homomorphism of Atiyah-Adams. More precisely they proved that if
$f$ is such a family then there is a homotopy invariant belonging to
the first (singular) cohomology group with $\Z_2$ coefficients whose
non-triviality insures the existence of a bifurcation point from the
trivial branch. This invariant is nothing but  the first
Stiefel-Whitney class of the analytical index and its non-triviality
can be geometrically interpreted as the non-orientability of the
index bundle.

We point out that all the quoted papers are non variational in
nature and as we shall see the situation in the variational case is
deeply different. For this reason, in this paper we restrict
ourselves to the case of families of self-adjoint Fredholm maps
which naturally arise in calculus of variation, by showing that in
this situation a more refined invariant arises. In this case in
fact, we will show that the bifurcation is related to the non
triviality of the {\em first Chern class of the index bundle\/}
instead to the first Stiefel-Whitney class. Moreover, by using this
invariant, we are able to prove a result about the {\em Lebesgue
covering dimension of the bifurcation set\/}.

As an application of our main result we give a quantitative result
on the  bifurcation set for a several parameter family of perturbed
geodesics on semi-Riemannian manifolds. We recall that perturbed
geodesics ($p$-geodesics) are trajectories  of a simple mechanical
system on a semi Riemannian manifold $M$ of finite dimension $n$.
They are critical points of the energy functional defined on the
space of paths on $M$ by
\begin{equation}\label{energia0}
E ( \gamma) \, = \, \int_{0}^{1} \frac{1}{2}g(\gamma'(x),
\gamma'(x))dx \, - \, \int_{0}^{1} V(x,\gamma(x)) dx
\end{equation}
\noindent  where $g$ is  a semi Riemannian metric on $M$  and $V$ is
a  time-dependent potential.  When the metric is indefinite, the
functional $E$  is of   strongly indefinite type and the ordinary
Morse index of its critical points is  infinite. In \cite{MPP}, we
associated to each non degenerate perturbed geodesic $\ga$   a
generalized Morse index $ \ispec (\ga) \in \Z$ and shown that it
coincides with  the conjugate index $\icon(\ga),$ an integer that
counts algebraically the total number of conjugate points along
$\gamma$.

We remark that very recently several kind  of bifurcation phenomena
have been studied in differential geometry and in  mechanics. (See
\cite{MPP}, \cite{MusPejPor07} \cite{PPT} and references therein).
In the Riemannian case, the change of qualitative properties of
closed geodesics under one parameter variation of the metric was
discussed for instance in \cite[Section 3.4]{KL1}.

Now in order to briefly describe our result we recall that given a
smooth branch $\sigma \equiv \{\sigma_\lambda \colon \lambda \in X\}
$ of perturbed geodesics, a {\em bifurcation point from the trivial
branch} $\sigma$  is a point $\la_*$ in $X$ such that arbitrary
close to the pair $(\la_*,\sigma_{\la_*}) $ there are pairs $(
\la,\ga_\la) $ where $\ga_\la$ is a $p$-geodesic with the same end
points as $\sigma_\la$ but not belonging to the branch $\sigma$.
Denoting by $c_1(\Ind_s(L_\gamma^\C))$ the first Chern class of the
analytical index associated to the family of hessians of the
$p$-geodesic functional at the critical section $\sigma$ (See
Sections \ref{sec:indexbundleechern} and \ref{sec:final} for a
precise definition), we proved the following.
\begin{thm}\label{mainbifgeointro} Let $ M$ be a smooth manifold and
$(g_\lambda,V_\lambda)_{\lambda\in X} $ be a family of
semi-Riemannain metrics and smooth time dependent potentials on $M$
smoothly parameterized by the compact connected orientable
$n$-dimensional manifold $X$. Let $ \sigma $ be a trivial branch of
$p$-geodesics of $ (M, g_\la, V_\la)_{\la \in X}$ such that for some
$\lambda_0 \in X$ the $p$-geodesic  $\sigma_{\lambda_0}$ is
non-degenerate (as critical point of the energy functional). If
\[
c_1(\Ind_s(L_\gamma^\C))\neq 0
\]
then the Lebesgue covering dimension of the bifurcation set is at
least $n-1$. Moreover either $B$ disconnect $X$ or it is
contractible to a point.
\end{thm}
We close this section by pointing out that by using our main results
similar quantitative results could be deduced also for Dirac
operators parameterized by  moduli spaces of connections. In fact in
this case it is known that the index bundle carries nontrivial
characteristic  Chern classes.


\section{Spectral flow and Index bundle}\label{sec:indexbundleechern}

The aim of this section is to briefly recall some  well-known
definitions and results  about the spectral flow for paths of
Fredholm quadratic forms. Furthermore we will construct the
analytical index for a family of selfadjoint Fredholm operators
smoothly parameterized by a compact topological space. In the final
part of the section we recall a result which relates the spectral
flow and the first Chern class of the analytical index of the family
which represents the key ingredient in order to prove our result.
Our main references are \cite{AtiSin69,BooWoj,
FitPej91,FitPejRec99,MusPejPor07,Pej01}.


\subsection{Spectral flow of Fredholm forms.}
\begin{defn}\label{def:quadraticaFredholm}
A {\em Fredholm quadratic form\/} on the real separable Hilbert
space $(\mathscr H, \langle \cdot, \cdot \rangle_\mathscr H)$ is a
function $q\colon \mathscr H\to \R$ such that there exists a bounded
symmetric form $b:= b_{q}\colon \mathscr H\times \mathscr H\to\R$
with $q(u)=b(u,u)$ and with $ \ker b $ of finite dimension.
\end{defn}
The space $\mathscr Q(\mathscr H)$ of all bounded quadratic forms is
a Banach space with the sup norm and the set $\mathscr Q_F(\mathscr
H)$ of all Fredholm quadratic forms is an open subset of $\mathscr
Q(\mathscr H)$ that is stable under perturbations by weakly
continuous quadratic forms.
\begin{defn}\label{def:d4}
Let  $q\colon (S^1, s_0)\to (\mathscr Q_F (\mathscr H), q_0)$ be a
based loop of Fredholm quadratic forms. We define the {\em spectral
flow of $q$ \/} as follows:
\[
\spfl(q,S^1):=\spfl(A_{q(t)},S^1),
\]
where $A_{q(t)}$ is the unique bounded linear self-adjoint Fredholm
operator such that
\[
\left\langle A_{q(t)}u,u\right \rangle_\mathscr H ={q(t)}(u) \qquad
\textrm{for all}\ \  u\in \mathscr H.
\]
\end{defn}
Given any differentiable loop of Fredholm quadratic forms $q\colon
S^1\rightarrow \mathscr Q_F(\mathscr H)$ then, for all $t \in S^1$,
the derivative $\dot{q}(t)$ is also a quadratic form. We say that a
point $t$ is a {\em crossing point\/} if $\ker b_{q(t)}\neq \{0\} $
and we say that the crossing point $t$ is {\em regular\/} if the
{\em crossing form\/} $\Gamma(q,t),$ defined as the restriction of
the derivative $\dot q(t)$ to the subspace $\ker b_{q(t)}$, is
non-degenerate. It is easy to see that regular crossings  are
isolated and that the property of having only regular crossing forms
is generic for paths in $\mathscr Q_F(\mathscr H)$. From theorem 4.1
in \cite{FitPejRec99} we deduce the following.
\begin{prop}\label{crossform}
If all crossing points of the path are regular then they are finite
in number and
\begin{equation}\label{crformula1}
\spfl(q,S^1)= \sum_i  \sgn \, \Gamma(q,t_i).
\end{equation}
\end{prop}
Based on the definition above, we are able to define the spectral
flow for a family of Fredholm quadratic forms defined on a separable
Hilbert bundle which is more suitable for the applications to
geometric variational problems.
\begin{defn}
A {\em generalized family of Fredholm quadratic forms\/}
parameterized by $X$ is a smooth function $ q\colon \widetilde
{\mathscr H}\to \R ,$ where $\widetilde {\mathscr H}$ is a real
separable Hilbert bundle over $X$ and $q$ is such that its
restriction $q_x$ to the fiber $\widetilde {\mathscr H}_x$ over $x$
is a (bounded) Fredholm quadratic  form.
\end{defn}
If $(X, x_0)=(S^1, s_0)$ and $q_{s_0} \in \Gl(\widetilde{\mathscr
H}_{s_0})$, we define the spectral flow $\spfl(q) := \spfl (q,S^1)$
of such a family $q$, by choosing a trivialization
\[
\Phi \colon S^1 \times \widetilde{\mathscr H}_{s_0} \to \widetilde
{\mathscr H}\colon (t, u) \mapsto \Phi_t(u)
\]
and by setting
\begin{equation} \label{sflow2}
\spfl(q) = \spfl (\widetilde q,S^1), \qquad \textrm{where}\ \
\widetilde q (t)[u] := q_t (\Phi_tu).
\end{equation}
From the cogredience property of the spectral flow the right hand
side of \eqref{sflow2} is independent of the choice of the
trivialization.
\begin{rem}\label{thm:spectralflowrealecomplesso}
We observe that the spectral flow could be defined equally well for
Fredholm Hermitian form. Furthermore we observe that  if $\mathscr
H^\C$ denotes the complexification $\mathscr H \otimes \C$ of the
real separable Hilbert space $\mathscr H$ then  for any path of
self-adjoint Fredholm operators $L:(S^1,s_0)\to \Fred_S(\mathscr H)$
we can define  its complexification $L^\C$. Now, as a direct
consequence of proposition \ref{crossform}, it is possible to show
that
\[
\spfl(L^\C, S^1)=  \spfl (L, S^1).
\]

\end{rem}

\subsection{Index bundle}

Given a compact space $X$ we will denote with $\KK(X)$ the
Grothendieck group of complex virtual vector bundle over $X$. By
construction $\KK(X)$ is the group completion of the abelian
semigroup $\Vect(X)$ of all complex vector bundle over $X$. Every
element in $\KK(X)$ can be written as a difference $[A]-[B]$, where
$[A]$ denotes the class of $A \in \Vect(X)$. Following
\cite{AtiSin69} for any continuous path  $L \colon X \to
\Fred(\mathscr H)$ where $\Fred(\mathscr H)$ denotes the set of all
complex linear bounded Fredholm operators we can associate a
homotopy invariant called the {\em index bundle\/} as follows.  By
compactness of $X$, there exists a finite dimensional subspace $V$
of $\mathscr H$ transverse to the family $L$ meaning that
\[
\Imm L_x +  V = \mathscr H, \qquad \forall \  x \in X.
\]
Thus the set
\[
Y=\{(x,v )\in X\times \mathscr H:\ \ \ L_x\,v \in V\}
\]
is the total space of a vector bundle over $X$ with fiber $Y_x :=
L^{-1}_x(V)$.
\begin{defn}
We define the {\em index bundle \/} of the family the element
$\Ind(L)\in \KK(X)$ which represents the stable equivalence class of
the  vector bundle $Y$ over $X$, i.e.:
\[
\Ind(L)= [Y]-[\Theta(V)],
\]
where $[\Theta(V)]$ is the trivial bundle over $X$ with fiber $V$.
\end{defn}
\begin{rem}
We observe that it  does not depend on the choice of the subspace
$V$ but just on the homotopy class of the family. (See for instance
\cite[Theorem 2.1]{ZKKP}).
\end{rem}
The index bundle of a family possesses the same homotopy addition
and logarithmic properties as the ordinary numerical index of a
Fredholm operator.  Denoting by $[X_1,X_2]$ the set of all free
homotopy classes of maps between the topological spaces $X_1$ and
$X_2$, a fundamental result proved by Atiyah and J\"anich shows that
the natural transformation $\Ind:[-,\Fred(\mathscr H)]\to \KK(-)$ is
an isomorphism.

Recall that $\KK$ can be extended to  two periodic generalized
cohomology theory by setting
\[
\KK^0(X)=\KK(X),\qquad  \textrm{and}\qquad \KK^{-1}(X)=\KK(SX).
\]
We will also need the following result.
\begin{thm}(\cite[Theorem B]{AtiSin69})
\begin{enumerate}
\item  From the topological point of view the
space $\Fred_S(\mathscr H)$ has three connected components denoted
by $\Fred^+_S(\mathscr H)$, $\Fred^-_S(\mathscr H)$ and
$\Fred^i_S(\mathscr H)$ respectively of all {\em essentially
positive, negative and indefinite,\/} where, essentially positive
(resp. negative) means that the spectrum is positive (resp.
negative) on some invariant subspace of $\mathscr H$ of finite
codimension.
\item The components $\Fred^+_S(\mathscr H)$, $\Fred^-_S(\mathscr
H)$ both are contractible.
\item Define a map
\[
\beta \colon \Fred_S^i(\mathscr H) \to \Omega \Fred (\mathscr H)
\]
where $\Omega \Fred(\mathscr H)$ denotes the loop of the space of
all Fredholm operators obtained by assigning to each $A \in
\Fred_S^i(\mathscr H)$ the path from $\Id$ to $-\Id$ in
$\Fred(\mathscr H)$ given by
\[
\beta(A):=  \Id \cos(\pi t)+ i A \sin(\pi t), \qquad  t \in [0,1]
\]
and closed to a loop by the standard continuation
\[
\Id\left(\cos(\pi t) +i \sin (\pi t)\right), \qquad t \in ]1,2].
\]
Then $\alpha$ is a homotopy equivalence.
\end{enumerate}
\end{thm}
From the above theorem we have:
\[
\pi[SX, \Fred] \simeq \pi[X, \Omega \Fred], \qquad
\textrm{and}\qquad \pi[X, \Fred^i_S(\mathscr H)] \simeq \KK^{-1}(X)
\]
and hence $\Fred_S^i(\mathscr H)$ is a classifying space for the
functor $K^{-1}$.
\begin{equation}\label{eq:omotopie}
\KK^{-1}(X):= \KK^0(SX) \simeq \pi[SX, \Fred(\mathscr H)] \simeq
\pi[X, \Omega \Fred(\mathscr H)].
\end{equation}
Now we are able to introduce the following key definition.
\begin{defn}\label{def:analiticalindexselfadjointcase}
Let $L \colon X \to \Fred^i_S(\mathscr H)$ be a continuous family of
 bounded linear self-adjoint Fredholm operators
parameterized by $X$.
The {\em analytical index\/} of the family is
the homotopy class of this map as an element of the group
$\KK^{-1}(X)$ namely:
\[
\Ind_s L := \Ind(*\,\beta\, L]
\]
where $*:\{X,\Omega Y\}\to \{SX,Y\}$ is defined by
$*(f)([x,t]):=f(x)(t)$.
\end{defn}
Denoting by $\bar c_1$ the first Chern class and by taking into
account its funtorial properties,
it is readily seen that it depends only on the stable equivalence
class of complex vector bundles and hence it induces a natural
transformation of functors
\begin{equation}\label{eq:Cherntrasfnaturale}
\bar c_1 \colon \KK^0(-) \longrightarrow H^2(-, \Z).
\end{equation}
We denote by $\sigma^K$ and $\sigma^H$ the suspension isomorphisms
in $K$-theory and in singular cohomology respectively and let us
consider the commutative diagram below. We define the natural
transformation $c_1: \KK^{-1}(X) \to H^1(X)$ by means of the
following commutative diagram:

\begin{equation*}
\xymatrixcolsep{2pc}\xymatrix{ K^0(SX) \ar[r]^{\bar c_1} &H^2(SX,\Z)
\\
 K^{-1}(X)\ar[u]^{\sigma^K} \ar[r] &
H^1(X,\Z)\ar[u]^{\sigma^H}}.
\end{equation*}
The following result, whose proof can be found in \cite{BooWoj},
relates the spectral flow with the analytical index.
\begin{prop}(\cite[Corollary 1.13]{BooWoj})\label{thm:Wojciekowski}
Let $L \colon S^1 \to \Fred^i_S(\mathscr
H)$ be a continuous path. Then we have:
\[
\spfl(L, S^1) =  c_1\,(\Ind_\alpha L)[S^1]
\]
where as usually $[\cdot ]$ denotes the fundamental class.
\end{prop}


\section{The main result}\label{sec:variationalframe}

Let $H$  be a real Hilbert space and let $X$ be an orientable
compact connected $n$-dimensional manifold. Let $f: X \times H \to
\R$ be a smooth family of functionals parameterized by $X$ and such
that $0$ is critical point of $f(x,\cdot)$ for any $x \in X$. We
will study bifurcation of critical points of $f$ from the trivial
branch $X \times \{0\}$.
\begin{defn}\label{def:biffamily}
A point $x_* \in X$ is called a bifurcation point of $\nabla f$ from
the trivial branch $X \times \{0\}$ is each neighborhood of
$(x_*,0)$ contains solutions of $\nabla f(x, h)=0$ with $h \neq 0$.
\end{defn}
\begin{notation}
We shall denote by $B:=\Bif(\nabla f)$ the set of all bifurcation
points of $\nabla f$ from the trivial branch.
\end{notation}
The map $x \mapsto d^2 f_x(0)$ defines a smooth family of bilinear
forms (with respect to the norm topology of all linear bounded
operators $\mathscr L(H)$); this  by Riesz representation theorem we
have that
\[
d^2f_x(0)[u,v]=\langle A_xu,v\rangle.
\]
Clearly $A:X \ni x \mapsto A_x\in \mathscr F(\mathscr H)$ is a
smooth map into the space of all bounded selfadjoint operators on
$H$. We will assume that
\begin{enumerate}
\item[(S)] for each $x \in X$ the operator $A_x$ is an indefinite Fredholm
operator smoothly dependent on $x \in X$.
\end{enumerate}
Thus it remains well-defined the smooth map $A: X \to \Fred^i(H): x
\mapsto A_x$.
\begin{mainthm}\label{thm:generalbif}
Let $X$ be a compact connected orientable $n$-dimensional manifold,
$H$ be a real separable Hilbert space and
\[
f: X \times H \to \R
\]
be a smooth family of of functionals parameterized by $X$ such that
$(S)$ holds. Assume that
\begin{enumerate}
\item[(i)] $c_1(\Ind_s A^\C)\neq 0$;
\item[(ii)] there exists some point $x_0 \in X\backslash B$ such
that $A_{x_0}$ is invertible.
\end{enumerate}
Then the Lebesgue covering dimension of $B$ is at least $n-1$.
Moreover either $B$ disconnect $X$ or it is contractible to a point.
\end{mainthm}
\proof First of all we observe that $B$ disconnects $X$ then the
dimension of $B$ is at least $n-1$ since the complement of a closed
subset of dimension less than $n-1$ is always connected. Therefore
let us assume that $X\backslash B$ is connected.

Let $A^\C$ be the complexified family. By assumption (i), there
exists a homology class $\gamma \in H^1(X, \Z)\neq 0$ such that
$\langle c_1(\Ind_s A^\C), \gamma\rangle \neq 0$. To prove the
statement we will show that the restriction to $B$ of the Poincar\'e
dual of $\gamma$ gives arise a nontrivial class in $\overline H^{(n
-1)}(B, \Z)$, i.e. $\iota^*(PD(\gamma))\neq 0$ where $\iota^*$ is
the homomorphism induced in cohomology by the inclusion $\iota: B
\hookrightarrow X$ and $\overline H^*(B, \Z):= \dirlim\, H^*(U, \Z)$
where $U$ ranges over all neighborhoods of $B$ in $X$. To do so we
consider the commutative diagram:
\begin{equation*}
\xymatrixcolsep{2pc} \xymatrix{ &  & H^{n-1}(X, \Z) \ar[r]^{\iota^*}
&\overline H^{n-1}(B,\Z)
\\
\dots \ar[r] & H_1(X\setminus B; \Z)  \ar[r]^{j_*}& H_1(X, \Z)
\ar[r]^{p_*} \ar[u] & H_1(X,X\setminus B,\Z) \ar[r]\ar[u]&
H_0(X\setminus B, \Z).}
\end{equation*}
By commutativity of the diagram above, and by taking into account
the exactness of the pair $(X, X \setminus B)$ in homology, the
restriction $\iota^*(PD(\gamma))$ of $PD(\gamma) $ to $B$ is dual to
$p_*(\gamma)$, where $p_* \colon H_1(X, \Z) \to H_1(X, X\backslash
B, \Z)$. By contradiction we assume that $p_*(\gamma)=0$; thus
$p_*(\gamma)=j_*(\beta)$, for some $\beta \in H_1(X\backslash B,
\Z)$ and for $j \colon X\backslash B \hookrightarrow X$. Since
$X\backslash B$ is connected, the Hurewicz homomorphism $\Pi \colon
\pi_1(X \backslash B, x_0) \to H_1(X \backslash B, \Z)$ is
surjective. Thus $\beta =g_*([S^1])$ for $g \colon (S^1, s_0) \to
(X\backslash B, x_0)$. Let $\bar g: (I, \partial I)\to (X, x_0)$ be
defined as the composition $j \circ g \circ \pi$ where:
\[
\pi: (I, \partial I) \to (S^1, s_0), \quad g:(S^1, s_0)\to
(X\backslash B, x_0),\quad  j:(X\backslash B, x_0)\to (X, x_0).
\]
Consider now the family of functionals $\bar f: I \times H \to \R$
defined by
\[
\bar f(\lambda, u) :=f(\bar g(\lambda), u).
\]
Then $d^2 \bar f_\lambda(0)= d^2 f_{\bar g(\lambda)}(0)$ and
therefore we denote with $\bar A_\lambda$  the corresponding family
of Hessian operators $\bar A:= A \circ \bar g$. Since $\lambda \in
I$ is a bifurcation point for $\nabla \bar f$ if and only if
$g(\lambda) \in B$ it follows that $\nabla \bar f$ is free of
bifurcation points. Moreover being $\bar A_0, \bar A_1$ invertible
operators and by taking into account \cite[Theorem 1]{FitPejRec99},
we can conclude that $ \spfl (\bar A, I)=0.$  But
\[
\spfl(\bar A, I)\ =\ \spfl (\bar A^\C, I)\ =  \spfl(A^\C \circ \bar
g, I)\ =\  \spfl (A^\C\circ jg, S^1),
\]
where the first equality is a direct consequence of remark
\ref{thm:spectralflowrealecomplesso}. By proposition
\ref{thm:Wojciekowski}, we have:
\[
\spfl(A^\C\circ j g, S^1)=\langle c_1(\Ind_s(A^\C\circ j g)) ,
[S^1]\rangle= \langle c_1(\Ind_s(A^\C)) ,j_*g_*( [S^1])\rangle =0.
\]
However by hypothesis $\langle c_1(\Ind_s(A^\C), \gamma\rangle \neq
0$ which is a contradiction. Therefore $p_*(\gamma)\neq0$ and hence
$\iota^*(PD(\gamma))$ restrict to a nontrivial class in $\overline
H^{n-1}(B, \Z)$. Now the conclusion follows by taking into account
the homological characterization of the Lebesgue covering dimension
(see for instance \cite[Chapter VIII, Theorem VIII.4]{HurWal48}) of
a topological space. \finedim

\section{A geometric variational framework: the perturbed geodesics case}\label{sec:final}

The aim of this section is to deduce some quantitative results for
the bifurcation set of a several parameter family of perturbed
semi-Riemannian geodesics. However in order to do so we need to
reformulate the previous result on Hilbert bundle. We start with the
following definition.
\begin{defn}
A {\em smooth family of Hilbert manifolds\/} $\{E_x \}_{x\in X}$
parameterized by $X$ and modeled on the real separable Hilbert space
$E$ is a family of manifolds of the form $ E_x =p^{-1}(x)$ where
$p\colon \widetilde E \to X$ is a smooth submersion of a Hilbert
manifold $\widetilde E$ onto $X$.
\end{defn}
By the implicit function theorem each fiber $E_x$ of the submersion
is a submanifold of $X$ of finite codimension (more precisely the
codimension is equal to the dimension of the compact manifold $X$).
For each $e\in E_x$, the tangent space $T_e E_x$ coincides with
$\ker \dd p_e$. Being $p$ a submersion the family of Hilbert vector
spaces $\TF(p) =\{\ker \dd p_x \colon x\in X\}$ is a Hilbert
subbundle of the tangent bundle $T\widetilde E$. $\TF(p)$ is the
bundle of tangents
along the fibers or the vertical bundle of the submersion $p.$\\
\noindent
 Thus given  a smooth functional $f\colon
\widetilde E \to \R$, it defines by restriction to the fibers of $p$
a {\em smooth family of functionals} $f_x \colon E_x \to \R $. We
will assume that there exists a smooth section $\sigma \colon X\to
\widetilde E$ of $p$ such that $\sigma (x)$ is a critical point of
the restriction $f_x$ of the functional $f$ to the fiber $E_x$ and
we will refer to $\sigma$ as the {\em trivial branch\/} of critical
points of the family $\{f_x : x\in X\}$.
\begin{defn}\label{def:ptodibifpiupara}
A point $x_* \in X$ is a {\em bifurcation point \/} of the family
$\{ f_x : x\in X\}$ from the trivial branch $\sigma(X)$ if there
exists a sequence $x_n \to x_*$ and a sequence $e_n \to \sigma(x_*)$
such that $p(e_n) =x_n$ and each $e_n$ is a critical point of $
f_{x_n}$ not belonging to $ \sigma (X)$.
\end{defn}
In what follows we shall denote by $h_x$ the Hessian of $f_x$ at the
point $\sigma(x)$. Our next assumption is the following
\begin{enumerate}
\item[\bf{(H)}] for each $x \in X,$ the Hessian $h_x$ is a strongly indefinite
self-adjoint Fredholm quadratic form.
\end{enumerate}
The family of Hessians $h_x$, for $x \in X$, defines a smooth
function $h$ on the total space of the pull-back bundle $\sigma^*
\TF(p)$ of the vertical bundle $\TF(p)$ by the map $\sigma\colon
X\to \widetilde E$ such that the restriction of $h$ to each fiber is
a Fredholm quadratic form. The function $ h\colon \sigma^*\TF(p) \to
\R$ is a {\em generalized family of Fredholm quadratic forms.\/}
According to the previous notation, we shall denote by $ h^\C$ the
complexification of the family $ h$ to the complexified Hilbert
subbundle $\sigma^*\TF(p)^\C= \sigma^*\TF(p)\otimes \C$. The
following result holds.

By  the vector bundle neighborhood theorem Theorem \cite[Appendix
A]{MusPejPor07},  there exist a trivial Hilbert bundle $\mathscr E
:= X\times E$ over $X$ and a fiber preserving smooth map $\psi\colon
\mathscr E \to \sigma^*\TF(p)$ such that $\psi(x ,0) = \sigma(x)$
and such that $\psi$ is a diffeomorphism of $\mathscr E$ with an
open neighborhood $\mathscr O $ of $\sigma(X)$ in $\sigma^*\TF(p)$.

Let $\widetilde f \colon X\times E  \to \R$ be defined by
$\widetilde f = f \circ \psi.$ Then $\widetilde f$ is a family of
smooth functionals on the real separable Hilbert space $E$. Since
the restriction $\psi_x$ of $\psi$  to the fiber   is a
diffeomorphism, we have that $ u \in E$ is a critical point of
$\widetilde f_x = f_x \circ \psi_x $ if and only if $  \psi_x(u)$ is
a critical point of $f_x$. In particular $0$ is a critical point of
$\widetilde f_x$ for each $x \in X.$ The Hessian $\widetilde h_x $
of $ \widetilde f_x$  at $0$ is given by $\widetilde h_x (\xi) = h_x
(d\psi_x(0))[\xi].$ Denoting by $\widetilde H_x (\xi, \eta)$ the
bilinear form associated to $\widetilde h_x$, by using the Riesz
representation theorem we have:
\[
\widetilde H_x (\xi, \eta)= \langle L_x \xi, \eta\rangle
\]
Clearly the map $x \mapsto L_x$ is a smooth map into the space of
all bounded selfadjoint operators on $H$; moreover assumption (H) is
equivalent to (S). Summing up the previous discussion and theorem
\ref{thm:generalbif}, we have the following.
\begin{mainthm}\label{thm:generalbif2}
Let $X$ be a compact connected orientable $n$-dimensional manifold,
$E$ be a real separable Hilbert space,  $p \colon \widetilde E \to
X$ a smooth submersion of a real separable Hilbert manifold
$\widetilde E$ over $X$ modeled on $E$ and $f \colon \widetilde E
\to \R$ be a smooth function. Assume that:
\begin{enumerate}
\item[(i)] $c_1(\Ind_s L^\C)\neq 0$;
\item[(ii)] there exists some point $x_0 \in X\backslash B$ such
that $L_{x_0}$ is invertible.
\end{enumerate}
Then the Lebesgue covering dimension of $B$ is at least $n-1$.
Moreover either $B$ disconnect $X$ or it is contractible to a point.
\end{mainthm}

\subsection{A geometric framework}

Let $M$ be a smooth connected  manifold of dimension $n$ and
$I:=[0,1]$. Let $g$ be a  semi Riemannian metric on $M$ and let $V
\colon I\times M\to \RR$ be  a smooth time-dependent potential. We
will denote by $D$ the associated Levi-Civita connection and by
$\Ddx$ the covariant derivative of a vector field along a smooth
curve $\gamma.$
\begin{defn}
A {\em perturbed geodesic} ($p$-geodesic) is a  solution $\gamma
\colon [0,1] \to M$ of the  second order differential equation
\begin{equation}\label{equ:pgeodesics}
\Ddx \, {\gamma '}(x)\, + \, \nabla\, V(x, \gamma(x)) \, = \, 0
\end{equation}
where  $\nabla$ is the gradient of a function on $M$ with respect to
the metric $ g.$
\end{defn}
Clearly if $V\equiv 0$ the above equation reduces to the classical
geodesic equation. Let $\Omega$ be the manifold of all $H^1$-paths
in $M$. It is well known that $\Omega$ is a smooth Hilbert manifold
modelled by $H^1([0,1];\RR^n).$  The tangent space $T_{\gamma }
\Omega$ at the point $\gamma$ can be identified in  a natural way
with the Hilbert space $H^1(\gamma) = \{ \xi \in H^1([0,1];TM)
\colon
      \tau \xi = \gamma \}$ of all $H^1$-vector fields along $\gamma$.
Here  $\tau  \colon TM \to M $ is the projection of  the tangent
bundle of $M$ to its base.

Any choice of a Riemannian (positive definite) metric on $M$ endows
$\Omega$ with an associated Riemannian  structure and hence with a
distance which makes $\Omega $ a metric space. The end-point map
\begin{equation*}
\pi \colon \Omega \to M\times M ; \, \quad \pi (\gamma ) = (\gamma
(0), \gamma(1))
\end{equation*} is  known  to be  a submersion and therefore for each
$(p,q) \in M\times M$ the fiber
\begin{equation*}
           \Omega_{p,q} = \{ \gamma \in \Omega \, :  \,
\gamma (0)=p , \, \gamma (1) = q \} \label{s2}
\end{equation*}
is a sub-manifold of codimension $2n$ whose tangent space $T_{\gamma
} \Omega_{p,q}$ is the subspace $H_0^1 (\gamma) $  of $H^1 (\gamma)
$ given by
\begin{equation*}
H^1_0 (\gamma ) = \{ \xi \in H^1 (\gamma ) \, : \, \xi (0)= \xi (1)
=0 \}.
\end{equation*}
Since  $\pi$ is a submersion,  the family of Hilbert spaces $H^1_0
(\gamma )$ is a Hilbert bundle $TF(\pi) = \ker \, T\pi $ over
$\Omega$, called the {\em bundle of tangents along the fibers}. The
reference for all the above,  and for everything else in this
section is  \cite{MPP}. Associated to each pair  $(g,V)$  there is
an {\em energy functional} $ \bar{E }\colon \Omega \to \RR$ defined
by
\begin{equation}\label{energia}
\bar E ( \gamma) \, = \, \int_{0}^{1} \frac{1}{2}g(\gamma'(x),
\gamma'(x))dx \, - \, \int_{0}^{1} V(x,\gamma(x)) dx.
\end{equation} The critical points of the restriction $ E= E_{p,q}$
of $\bar E$ to the sub-manifold $\Omega_{p,q}$ are precisely the
$p$-geodesics through $p, q$ with kinetic energy  $\frac{1}{2}
g(\ga',\ga')$ and potential energy $V$. To be  more precise,
critical points  of the energy functional are weak solutions of
\eqref{equ:pgeodesics}, which turn out to be smooth classical
solutions  by elliptic regularity. In \cite{MPP}, we  associated two
integers to each  non degenerate critical point of the energy
functional. The first integer is the generalized Morse index
$\ispec(\gamma)$, defined as the negative of the spectral flow of
the Hessians of the energy functionals  along the path $\tilde\ga$
canonically induced by $\gamma$ on  $\Omega.$  Since several points
related to the construction of the generalized Morse index will be
used, we will briefly recall their definitions.

The linearization of the boundary value problem
\[
\left\{\begin{array}{ll}\Ddx \, {\gamma '}(x)\, + \,\nabla\,
    V(x, \gamma(x)) \, = \, 0 &  \\\ga(0)=p, \quad \ga(1)=q & \end{array}\right.
\]
    at the critical point $
\ga $ is the  equation of Jacobi fields

\begin{equation}\label{jaceq}
\Ddtt\xi (t) + R(\gamma ' (t) , \xi (t)) \gamma '(t) + D_{\xi
(t)}\nabla V(t, \gamma(t))=0,
\end{equation} subjected to Dirichlet boundary conditions
$\xi(0)=0= \xi(1).$ Here $\xi$ is a vector field along $\ga,$ $R$ is
the curvature tensor of
    the connection $D$ and $ D_{\xi(x)} \nabla V(x, \ga  (x))$ is the
Hessian of $V(x,-)$ with respect to the metric $g,$ viewed as a
symmetric  operator of $T_{\ga(x)} M.$ Its weak solutions
($H_0^1$-{\em Jacobi fields}) are vector fields $\xi \in
H^1_0(\ga)$ such that for any $\eta \in H^1_0(\ga)$
\begin{eqnarray} \label{Hessi} 0&=& H_\ga (\xi, \eta) =
\int_{0}^{1}  g( \Ddx \xi(x), \Ddx \eta(x)) \, dx +\\\nonumber
&-&\int_{0}^{1}g(R( \ga'(x), \xi(x)) \ga'(x) + D_{\xi(x)} \nabla
V(x, \ga (x)), \eta(x))  \, dx.
\end{eqnarray}
Thus,  Jacobi fields in $H^1_0(\ga)$ are the elements of the kernel
of the bilinear form $H_\ga$. The Hessian of $E_{p,q}$ at $\ga$ is
the  quadratic form $h_\ga \colon  H^1_0(\ga) \to \R$ associated to
$ H_\ga$. Since the embedding of  $H^1_0(\ga)$ into $L^2(\ga)$  is
compact, the form
\[ c(\xi)  = \int_{0}^{1}g(R( \ga'(x), \xi(x)) \ga'(x) +
D_{\xi(x)} \nabla V(x, \ga (x)), \xi(x))  \, dx\]  is weakly
semi-continuous. Being a weakly semi-continuous perturbation of a
non degenerate form,  $h_\ga $ is a bounded Fredholm quadratic form.
In particular  $\ker h_\ga  \equiv \{ \xi |\, H_\ga (\xi,\eta ) = 0
{ \,\rm\,  for \  all \,} \eta\}$ is finite dimensional. The form
$h_\gamma$ is non degenerate  if and only if the instant $1$ is not
conjugate to $0$ along $\gamma$.  If this holds we will say that the
$p$-geodesic is non degenerate. Each  $p$-geodesic $\ga$  induces in
a canonical way a path $ \tilde \ga $  on the manifold $\Omega $ of
all paths on $M$. The canonical path  $ \tilde \ga\colon  \to
\Omega$  associated to $\ga$ is defined by $\tilde\ga (t)(x) =
\gamma (t \cdot x )$, $x \in [0,1]$. Each $\tilde\ga_t \equiv
\tilde\ga(t) $ is a $p$-geodesic for $(g, V_t )$, where $V_t(x,m) =
\lambda^2V(t x,m)$. From the above discussion it follows that the
Hessian $h_t$ of $E_t \equiv E_{p, \ga(t)}$ at $\tilde\ga_t$ is
degenerate  if and only $\gamma (t)$ is a conjugate  point to $
p=\gamma (0)$ along $\gamma$. If $\ga$ is non degenerate  the path
$\tilde\ga$ has non degenerate end points. The family $h_t,\ t\in
[0,1]$ defines a function $h$ on the total space of the pull-back
bundle $ \tilde\ga^*( TF)$ over $[0,1]$ such that its restriction to
each fiber is a bounded quadratic form and such that $h_0, h_1$ are
non degenerate.  This is what we called in \cite{MPP} an  admissible
family of Fredholm quadratic forms. The spectral flow of such a
family was defined in \cite[Section 2]{MPP}.
\begin{defn}
We define the  {\em generalized Morse index}  of a $p$-geodesic
$\ga$ as:
\[
\ispec(\gamma)= -\spfl (h).
\]
\end{defn}
\begin{defn}
Let $X$ be a compact, orientable connected $n$-dimensional manifold
and let $p \colon \mathcal S_\nu ^{2}(M)  \to M$ be the bundle  of
all non degenerate symmetric two forms of index $\nu$ on the tangent
bundle $TM$. A {\em several parameter family of semi Riemannian
metrics on $M$} is a smooth map $g\colon X \times M \to \mathcal
S_\nu^{2}(M)$ such that $pg =\pi,$ where $ \pi $ denotes the
projection onto the second factor.
\end{defn}
Each   $ g_{\lambda}= g (\lambda, -)$ is then a semi Riemannian
metric on $M$ of index $\nu.$ Using the Koszul formula  it is easy
to  see that the  family  of associated Levi-Civita connections  $
D^{\lambda}$  is smooth with  respect to the $\lambda$ variable in
an obvious sense. If  $V\colon X \times [0,1]\times  M \to \RR$ is a
smooth function, we consider $V$ as a smooth one parameter family of
time-dependent {\em potentials} $V_\lambda\colon I \times M \to \RR
$ defined by $V_\lambda (x,u) = V(\lambda, x, u).$ The data $( M,
g_\lambda,V_\lambda)_{\lambda\in X} $ define a several parameter
family of time-dependent  mechanical systems.

We shall also assume the  existence of a several parameter family of
known $p$-geodesics. In other words, given a smooth family
$(g_\lambda,V_\lambda)_{\lambda\in X} $ as above,  we assume that
there exists  a smooth map $\sigma\colon X\times I \to M$ such that
for each $\la \in X$ ,  $ \sigma_\lambda (x)\equiv \sigma (\lambda ,
x)$ is a $p$-geodesic corresponding to the mechanical system
$(g_\lambda, V_\lambda)$.

We will refer to the family $(\sigma_\lambda )_{\lambda \in X}$ as
{\em the trivial branch } of perturbed geodesics.
\begin{defn}\label{def:punto di bif} A point $\lambda_* $ is called
a {\em bifurcation point} for $p$-geodesics  from the trivial branch
$\sigma$ if there exists a sequence $(\lambda_n, \gamma_n) \to
(\lambda_*,\sigma(\lambda_*))$ in $X\times \Omega$ such that
$\gamma_n$ is a weak solution of  equation \eqref{equ:pgeodesics}
with  boundary conditions
\begin{equation} \label{bc}
\gamma_n(0) = \sigma _{\lambda_n}(0) , \quad  \gamma_n(1) =
\sigma_{\lambda_n}(1)
\end{equation} such that $ \gamma_n$ does not belong to  $
\sigma (X)$.
\end{defn}
For each $\lambda \in X$, let  $E_\lambda$ be the restriction of the
energy functional associated to the data $(g_\lambda , V_\lambda)$
to the sub-manifold $\Omega_\la \equiv \Omega_{\sigma_\lambda (0),
\sigma_\lambda(1)}$. The trivial branch $\sigma$ can be viewed  as a
smooth path  $\sigma \colon I \to \Omega$ of critical points
$\sigma(\lambda) = \sigma_\la$ of $E_\la.$ We set
\[
c_1(\Ind_s h_\gamma^\C):=c_1(\Ind_s(L_\gamma^\C)).
\]
As a direct consequence of theorem \ref{thm:generalbif2} the
following result holds true.
\begin{thm}\label{mainbifgeo} Let $ (M, g_\lambda,V_\lambda)_{\lambda\in X} $ be as above
and let  $ \sigma $ be a trivial branch  of  $p$-geodesics  of $ (M,
g_\la, V_\la)_{\la \in X}$ such that for some $\lambda_0 \in X$,
$\sigma_{\lambda_0}$ is non-degenerate. If
$c_1(\Ind_s(L_\gamma^\C))\neq 0$ then the Lebesgue covering
dimension of the bifurcation set is at least $n-1$. Moreover either
$B$ disconnect $X$ or it is contractible to a point.
\end{thm}


\begin{thebibliography}{99}
\bibitem{AtiSin69} M.F Atiyah, I. M. Singer. {\em Index theory for
skew-adjoint Fredholm operators.\/} Publ. Math. de I.H.E.S. {\bf 37}
(1969), 5--26.

\bibitem{BooWoj} B. Booss, K. Wojciechowski. {\em Desuspension of splitting elliptic
symbols I.\/} Annals of Global Analysis and Geometry {\bf 3} (1985),
337--383.


\bibitem{FitPej91} P.M. Fitzpatrick - J. Pejsachowicz. {\em Nonorientability of the Index bundle and Several-Parameter
bifurcation},  Journal of Functional Analysis {\bf 98} (1991), N.1
42--58.


\bibitem{FitPejRec99} P.M. Fitzpatrick - J. Pejsachowicz - L. Recht.
{\em Spectral flow and bifurcation of critical points of
strongly-indefinite functional. Part I. General theory}, J.
Functional Analysis {\bf 162} (1999), 52--95.


\bibitem{HurWal48} W. Hurewicz, H. Wallman {\em Dimension theory.\/} Princeton,
Princeton University Press 1948.

\bibitem{KL1} W. Klingenberg.  Riemannian Geometry, {\em de Gruyter,
New York, (1995)}.


\bibitem{Kui65} Kuiper, Nicolaas H. {\em The homotopy type
of the unitary group of Hilbert space.\/} Topology 3 1965 19--30. 

\bibitem{MPP} M. Musso - J. Pejsachowicz - A. Portaluri. {\em A Morse
Index Theorem and bifurcation for perturbed geodesics on
Semi-Riemannian Manifolds}, Topol. Methods Nonlinear Anal. {\bf 25}
(2005), no. 1, 69--99.

\bibitem{MusPejPor07}M. Musso, J. Pejsachowicz, A. Portaluri
{\em Bifurcation of perturbed geodesics on semi-Riemannian
manifold.\/} Esaim-Cocv {\bf 13} (2007), n.3, 598--621.

\bibitem{Pej88} J. Pejsachowicz {\em $K$-theoretic methods in bifurcation theory.\/} Contemporary Mathematics
{\bf 72}, 1988.

\bibitem{Pej01} J. Pejsachowicz {\em Index Bundle. Leray-Schauder reduction and bifurcation of
solutions of nonlinear elliptic boundary value problems.\/}
Topological Methods in Nonlinear Analysis {\bf 18}, (2001),
243--267.

\bibitem{PPT} P. Piccione - A. Portaluri - D.V.Tausk. {\em Spectral
flow, Maslov index and bifurcation of semi- Riemannian geodesics},
to appear in Annals of Global Analysis and Geometry.




\bibitem{ZKKP} M. G. Zaidenberg, S.G. Krein, P. A. Kuchment, A. A. Pankov {\em Banach Bundles and linear
operators.\/} Russian Math. Surveys {\bf 30} (1975), 115-175.

\end{thebibliography}
\end{document}